\documentclass[a4paper,10pt,twoside]{article}

\usepackage[austrian,american]{babel}
\usepackage[latin1]{inputenc}
\usepackage[intlimits,leqno]{amsmath}
\usepackage{amssymb}
\usepackage{amsthm}
\usepackage{graphicx,color}
\usepackage{rotating}
\usepackage{fancyhdr}
\usepackage[nottoc]{tocbibind}
\usepackage{makeidx}
\usepackage{epsfig}

\makeindex

\usepackage{latexsym}

\usepackage{pstricks,pst-node,amsfonts}
\definecolor{darkgray}{rgb}{0.5,0.5,0.5}
\definecolor{lightgray}{rgb}{0.8,0.8,0.8}
\newcommand\R{{\mathbb{R}}}

\newcommand\C{{\mathbb{C}}}
\newcommand\N{{\mathbb{N}}}
\newcommand{\A}{\mathcal{A}}

\newcommand{\X}{\mathcal{X}}

\newcommand{\F}{\mathcal{F}}

\renewcommand{\H}{\mathcal{H}}

\newcommand{\x}{\mathbf{x}}

\newcommand{\bea}{\begin{eqnarray}}
\newcommand{\eea}{\end{eqnarray}}
\newcommand{\bean}{\begin{eqnarray*}}
\newcommand{\eean}{\end{eqnarray*}}

\newcommand{\be}{\begin{equation}}
\newcommand{\ee}{\end{equation}}

\newcommand{\ba}{\begin{array}}
\newcommand{\ea}{\end{array}}

\newtheorem{theo}{Theorem}[section]
\newtheorem{prop}[theo]{Proposition}
\newtheorem{coro}[theo]{Corollary}
\newtheorem{lemm}[theo]{Lemma}

\newtheorem{rema}[theo]{Remark}

\renewcommand{\d}{{\rm d}}

\makeatletter
\newcommand\reqref[2]{\textup{\tagform@{\ref{#1}--\ref{#2}}}}
\makeatother

\begin{document}

%\definecolor{dark}{gray}{.5}
%\definecolor{light}{gray}{.8}

%------------------------ title --------------------------

\title{Causality analysis of waves and wave equations obeying attenuation
%Modeling of wave equations
%obeying frequency dependent attenuation laws
%for thermoacoustic Tomography
}

\author{Richard Kowar\\
Department of Mathematics, University of Innsbruck,\\
Technikerstrasse 21a/2, A-6020 Innsbruck, Austria\\
richard.kowar@uibk.ac.at}
%\email{richard.kowar@uibk.ac.at}
%\affiliation{Department of Mathematics, University of Innsbruck,\\
%Technikerstrasse 21a/2, A-6020 Innsbruck, Austria
%}

\date{January 9, 2009}

\maketitle

\begin{abstract}
In this paper we show that the standard causality condition for attenuated waves, i.e. the Kramers-Kronig relation that relates the
attenuation law and the phase speed of the wave, is necessary but not sufficient for causality of a wave. By causality of a wave we understand
the property that its wave front speed is bounded. Although this condition is not new, the consequences for wave attenuation have not been
analysed sufficiently well. We derive the wave equation (for a homogeneous and isotropic medium) obeying attenuation and causality and with a
generalization of the Paley-Wiener-Schwartz Theorem (cf. Theorem 7.4.3. in~\cite{Ho03}), we perform a causality analysis of waves obeying the
frequency power attenuation law. Afterwards the causality behaviour of Szabo's wave equation (cf.~\cite{Szabo94}) and the thermo-viscous wave
equation are investigated. Finally, we present a generalization of the thermo-viscous wave equation that obeys causality and the frequency
power law (for powers in $(1,2]$ and) for small frequencies, which we propose for Thermoacoustic Tomography.
\end{abstract}

\textbf{AMS:} 45K05, 35Q72, 74J05, 42A38, 42A85

%------------------- begin introduction -----------------------------

%\cleardoublepage
%\pagenumbering{arabic}

\section{Introduction}
\label{sec-intro}

In physics an attenuated wave is modeled by replacing the real frequency dependent wave number by a complex frequency dependent wave number.
For an \emph{attenuated spherical wave} with origin $(\x,t)=(\mathbf{0},0)$, this yields~\cite{FeLeSa63,Ro07,KiFrCoSa00}
\begin{equation}\label{patt}
\begin{aligned}
   p_\alpha(\x,t)
      = \frac{1}{\sqrt{2\,\pi}} \int_\R \frac{e^{-\alpha(\omega)\,|\x|}
            \,e^{-i\,\omega\,\left(t-\frac{|x|}{c(\omega)}\right)}}{4\,\pi\,|\x|}\,\d \omega\,,
%      = \frac{g_\alpha(\x,t-\frac{|\x|}{v_0(|\x|)})}{4\,\pi\,|\x|}\,,
\end{aligned}
\end{equation}
where $\alpha=\alpha(\omega)>0$ is called the \emph{attenuation law}, which is assumed to be a positive real-valued even function, and
$c=c(\omega)$ is called the \emph{phase speed}. The complex wave number is given by $k(\omega)= \frac{\omega}{c(\omega)} + i\,\alpha(\omega)$.
With the notion
\begin{equation}\label{defalpha*}
   \alpha_*(\omega) := \alpha(\omega)-i\,\left(\frac{\omega}{c(\omega)}-\frac{\omega}{c_0}\right)
                     = -i\,k(\omega) + i\,\frac{\omega}{c_0}
  \qquad (c_0>0 \mbox{ const.})\,,
\end{equation}
the standard causality requirement for waves of the form~(\ref{patt}) read as follows (cf.~\cite{Wa00} and Section 3.3 in~\cite{BeWo66})
\begin{equation}\label{KKst0}
\begin{aligned}
   \mbox{Re}(\alpha_*(\omega)) &= -\mbox{Im}(\H(\alpha_*(\omega)))\,\quad\mbox{ and } \quad
   \mbox{Im}(\alpha_*(\omega)) &= \mbox{Re}(\H(\alpha_*(\omega)))\,.
\end{aligned}
\end{equation}
These relations are the \emph{Kramers-Kronig relations} for $\alpha_*$  and are equivalent to the relation (cf.~\cite{Wa00})
\begin{equation}\label{KKst}
\begin{aligned}
     \frac{\omega}{c(\omega)}-\frac{\omega}{c_0} = -\H(\alpha(\omega))\,.
\end{aligned}
\end{equation}
Here $\H$ denotes the \emph{Hilbert transform}.

In this paper we require that every wave has a finite front speed $v_\alpha$.  Of course, this is not a new requirement (cf.~\cite{Br60}), but
the consequences for wave attenuation have not been analysed sufficiently well. If a wave $p_\alpha$ of the form~(\ref{patt}) has a finite
front speed $v_\alpha$ bounded from above by a constant $v_B>0$, then the distribution
\begin{equation*}
\begin{aligned}
   g_\alpha(\x,t) :=  4\,\pi\,|\x|\,p_\alpha\left(\x,t+\frac{|\x|}{v_B}\right) \qquad \x\in\R^3,\,t\in\R\,
\end{aligned}
\end{equation*}
is causal for any $\x\in\R^3$, i.e. its support lies in $[0,\infty)$ for any fixed $\x\in\R^3$. Conversely, if $g_\alpha(\x,\cdot)$ is causal
for any $\x\in\R^3$, then the front speed of the wave $p_\alpha$ is bounded from above by $v_B>0$. Let $\hat g_\alpha(\x,\omega)$ denote the
\emph{Fourier transform}\footnote{ If $f\in \mathcal{S}(\R)$, then $\hat f(\omega) :=\F(f)(\omega):=\frac{1}{\sqrt{2\,\pi}}\,\int_\R
e^{-i\,\omega\,t} f(t)\d t$ and $\check f :=\F^{-1}(f)$. This convention implies $\widehat{f\,g}=\frac{1}{\sqrt{2\,\pi}}\,\hat f * \hat g$ and
$\hat f\,\hat g=\frac{1}{\sqrt{2\,\pi}}\,\widehat{f*g}$.} of $g_\alpha(\x,t)$. Then $g_\alpha$ is causal if and only if the Kramers-Kronig
relations for $\hat g_\alpha$ holds (cf. Section 3.3 in~\cite{BeWo66}), i.e.
%\footnote{these relations hold also for the inverse Fourier transform $\check g_\alpha$ of $g_\alpha$. }
\begin{equation}\label{KKg}
\begin{aligned}
   \mbox{Re}(\hat g_\alpha) = -\mbox{Im}(\H(\hat g_\alpha))\,\quad\mbox{ and } \quad
   \mbox{Im}(\hat g_\alpha) = \mbox{Re}(\H(\hat g_\alpha))\,,
\end{aligned}
\end{equation}
If $g_\alpha(\x,\cdot)$ is a Schwartz function then these conditions imply the Kramers-Kronig relations~(\ref{KKst0}) for $\alpha_*$.
This can be seen as follows: Let $\x\in\R^3$ be arbitrary but fixed and let $*_t$ denote the time-convolution. For a spherical
wave of the form~(\ref{patt}) with $c_0:=v_B$ we have
\begin{equation}\label{Fst}
\begin{aligned}
   g_\alpha(\x,t)
     &=   \frac{1}{\sqrt{2\,\pi}}\,\mathcal{F}\left( e^{-\alpha_*(\omega)\,|\x|}\right)(t)
     \qquad (\x\in\R^3,\,t\in\R)
\end{aligned}
\end{equation}
with $\alpha_*$ defined as in~(\ref{defalpha*}).
Since the causality of $g_\alpha(\x,\cdot)$ implies the causality of
\begin{equation}\label{nablagst}
\begin{aligned}
  \nabla g_\alpha(\x,t)
     =  -\frac{\x}{|\x|}\,\hat \alpha_*(t) *_t g_\alpha(\x,t)\,,
\end{aligned}
\end{equation}
$\hat\alpha_*(t)$ must be causal.  But the causality of $\hat\alpha_*(t)$ means that the Kramers-Kronig relations~(\ref{KKst0}) for
$\alpha_*(\omega)$ hold. We note that in general, the causality of $\hat\alpha_*$ or $\nabla g_\alpha$ does not imply the causality of
$g_\alpha$ (cf. Section~\ref{sec-causcheck}).

The phase speeds induced by the frequency power laws $\alpha(\omega)=\alpha_0\,|\omega|^\gamma$ with $\gamma\in (\R^+\backslash \N)\cup\{1\}$ 
and the 
Kramers-Kronig relation~(\ref{KKst}), are derived in~\cite{Wa00,WaMoMi00,WaHuMoMi03,WaHuMoBrMi00,Szabo95}. For example if $\gamma\in
(0,\infty)$ and $\gamma\not\in\N$, then the phase speed reads as follows
\begin{equation*}
\begin{aligned}
%    \alpha(\omega) &= \alpha_0\,|\omega|^\gamma\,,\\
   \frac{1}{c(\omega)} - \frac{1}{c(\omega_0)} &= \alpha_0\,\tan\left(\frac{\pi}{2}\,\gamma\right)\,
                                   (|\omega|^{\gamma-1} - |\omega_0|^{\gamma-1})
 \qquad (\omega\in\R)\,,
\end{aligned}
\end{equation*}
which implies
\begin{equation}\label{al}
\begin{aligned}
\alpha_*(\omega)
    = \frac{\alpha_0\,(-i\omega)^\gamma}{\cos(\frac{\pi}{2}\,\gamma)}
       + i\alpha_0\, \tan\left(\frac{\pi}{2}\,\gamma\right)\,|\omega_0|^{\gamma-1} \,\omega\,.
\end{aligned}
\end{equation}
In Section~\ref{sec-causcheck} it will be shown that the wave~(\ref{patt}) with this $\alpha_*$ is causal only if $\omega_0=0$ and $\gamma\in
(0,1)$. This shows that the Kramers-Kronig relation~(\ref{KKst}) for $\alpha_*$ is necessary but not sufficient. Moreover, this shows that the
frequency power law for $\gamma\geq 1$ is not an \emph{admissible} attenuation law.
%In Section~\ref{sec-gamgen} we derive a wave equation that obeys the power frequency law ($\gamma\in (1,2]$) for sufficiently
%small frequencies and admits causality.
%Indeed, we propose this equation for Thermoacoustic Tomography, since it is said that soft tissue has a $\gamma-$value
%$\leq 2$.

A standard calculation shows that the wave~(\ref{patt}) with $\alpha_*$ defined as in~(\ref{al}) satisfies the following wave
equation:
\begin{equation}\label{standwaveeq}
\begin{aligned}
  \nabla^2 p_\alpha
       &-\left[\frac{\alpha_0}{\cos(\frac{\pi}{2}\,\gamma)}\,D_t^\gamma
         + \frac{1}{c_0} \,\frac{\partial}{\partial t} \right]^2 \, p_\alpha
      = -\delta(\x)\delta(t) \,,
\end{aligned}
\end{equation}
where $\alpha_0>0$  and $D_t^\gamma$ denotes the \emph{Riemann-Liouville fractional derivative} defined by~(\cite{KiSrTr06,Po99}))
\begin{equation}\label{defDtga}
        \F^{-1}(D_t^\gamma(f))(\omega) := (-i\,\omega)^\gamma\, \check f(\omega)\,.
\end{equation}
We note that the causality of the kernel of $D_t^\gamma$ implies the Kramers-Kronig relation~(\ref{KKst}) and vice versa.

A different wave equation was derived by Szabo in~\cite{Szabo94}). A comparision of the dispersion relations of the $1d-$ thermo-viscous wave
equation and the $1d-$ electromagnetic wave equation for conducting media lead Szabo to the following dispersion relation $k(\omega)^2 =
\frac{\omega^2}{c_0^2} + i\,2\,\frac{\omega}{c_0}\,\alpha_0\,|\omega|^\gamma$ with $\gamma>0$, which implies the wave equation:
\begin{equation}\label{szaboseq}
\begin{aligned}
    \nabla^2 p_\alpha -\frac{1}{c_0^2}\frac{\partial^2 p_\alpha}{\partial t^2}
       + L_\gamma *_t p_\alpha
            = -\delta(\x)\delta(t)\,.
\end{aligned}
\end{equation}
For example, if $\gamma>0$ and $\gamma\not\in\N$, then
\begin{equation}\label{Lgamma}
\begin{aligned}
  L_\gamma = - \frac{2\,\alpha_0}{\sqrt{2\,\pi}\,\cos(\frac{\pi}{2}\,\gamma)\,c_0}\,
         \F\left((-i\,\omega)^{\gamma+1}\right)\,.
%= - \frac{1}{\sqrt{2\,\pi}}\,
%         \F\left(\frac{2\,\alpha_0\,(-i\,\omega)^{\gamma+1}}{\cos(\frac{\pi}{2}\,\gamma)\,c_0}\right)\,.
%\qquad \mbox{ if $\gamma>0$ and $\gamma\not\in\N$}\,.
\end{aligned}
\end{equation}
Since $L_\gamma$ is the kernel of the operator $\frac{2\,\alpha_0}{\cos(\frac{\pi}{2}\,\gamma)\,c_0}\,\frac{\partial }{\partial
t}\,D_t^\gamma$, Szabo's equation can be obtained by neglecting the term $\frac{\alpha_0^2}{\cos^2(\frac{\pi}{2}\,\gamma)}\,D_t^{2\,\gamma}$
in equation~(\ref{standwaveeq}). In Section~\ref{sec-szabo} we show that Szabo's equation admits causality only if $\gamma\in (0,1)$.
And only in this case Szabo's equation is a small frequency approximation of equation~(\ref{standwaveeq}). \\

\noindent This paper is organized as follows: The general properties of attenuated waves assumed in this paper are introduced and discussed in
Section~\ref{sec-genmodel}. With these assumptions we derive the general wave equation for a homogeneous and isotropic medium obeying
attenuation and causality in Section~\ref{sec-genwaveeq}. The first causality analysis of attenuated waves is performed in
Section~\ref{sec-causcheck} for the case of the frequeny power law.  Afterwards the causality analysis of Szabo' equation
(Section~\ref{sec-szabo}) and the thermo-viscous wave equation (Section~\ref{sec-caustv}) are perfomed. Finally, a generalization of the
thermo-viscous wave equation, is derived in Section~\ref{sec-gamgen}, which admits causality and obeys the power frequency law for
sufficiently small frequencies and $\gamma\in (1,2]$.

\section{General properties of attenuated waves}
%General model for attenuated waves in homogeneous and isotropic media satisfying causality}
%\section{Basic properties of attenuated waves in homogeneous and isotropic media satisfying causality}
\label{sec-genmodel}

In this section we postulate the basic properties of attenuated waves propagating in \emph{homogeneous isotropic media} and
infer the structure of attenuated waves. The main goal of this section is to clarify our assumptions on which we base
our casuality analysis of attenuated waves.

\subsection*{Assumptions}

Let $p_0$ denote the solution of the standard wave equation $\nabla^2 p_0 -\frac{1}{c_0^2}\,\frac{\partial^2 p_0}{\partial t^2}
     = -f$
%\begin{eqnarray}\label{pwaveeq2}
%    \nabla^2 p_0 -\frac{1}{c_0^2}\,\frac{\partial^2 p_0}{\partial t^2}
%     = -f\,
%\end{eqnarray}
with $p_0|_{t<0}=0$ and $\frac{\partial p_0}{\partial t}|_{t<0}=0$. Let $\A_\alpha$ denote the map that associates to each source term $f$ the
corresponding attenuated wave $p_\alpha$. $\alpha=0$ means no attenuation, i.e. $p_0=\A_0(f)$. We assume a \emph{homogeneous isotropic
medium}, which implies that
\begin{equation*}
\begin{aligned}
     \mbox{$\A_\alpha$ commutes with all translations in $\x$ and $t$}\,.
\end{aligned}
\end{equation*}
We also assume that for $\alpha\not=0$: %attenuation smoothes the wave shape, i.e.
\begin{equation*}
\begin{aligned}
   &\mbox{$\A_\alpha$ is a linear continuous mapping} \\
   &\mbox{that maps  $C_0^\infty(\R^4)$ into $\mathcal{S}(\R,C^\infty(\R^3))$}\,,
\end{aligned}
\end{equation*}
and that
\begin{equation*}
\begin{aligned}
    \A_\alpha(\delta(\x)\,\delta(t))\in \mathcal{S}(\R,C(\R^3))\,.
\end{aligned}
\end{equation*}
Here $\mathcal{S}$ denotes the space of Schwartz functions, $\delta(\x)$ denotes the $3D-$delta distribution and $\delta(t)$ denotes the
$1D-$delta distribution. The last two assumptions take into account that wave attenuation smoothes and decreases the wave.
%It is known that $\A_0$ is a linear and continuous mapping that maps
%$C_0^\infty(\R^3\times\R)$ into $C^\infty(\R,C^\infty(\R^3))$,

\subsection*{Superposition law}

Since the operator $\A_\alpha$ satisfies the assumptions of Theorem~4.2.1 in~\cite{Ho03}, $\A_\alpha$ is a space-time convolution
operator with kernel $G_\alpha(\x,t):=\A_\alpha(\delta(\x)\,\delta(t))$.  Since $\delta(\x)\,\delta(t)$ is the source term of a
spherical wave for $\alpha=0$, we interprete $G_\alpha$ as an attenuated spherical wave and thus
\begin{equation}\label{Aalpha}
\begin{aligned}
 %    p_{att}(\x,t)  =
           \A_\alpha(f)(\x,t)
           = G_\alpha(\x,t)*_{\x,t}f(\x,t)\,
%            = \int_{\R^3}\int_R  G_\alpha(\x-\x',t-t')\,f(\x',t')\d \x' \d t'
\end{aligned}
\end{equation}
is nothing else but the superposition law of attenuated waves. In analogy to partial differential equations and linear system
theory we call $G_\alpha$ the \emph{Green function} of the wave model. \\

\subsection*{Causality condition}

We require that any attenuated spherical wave $G_\alpha$ has a positive finite wave front speed $v_\alpha$, which is equivalent to the
requirement that
\begin{equation}\label{defg}
\begin{aligned}
   g_\alpha(\x,t) :=  4\,\pi\,|\x|\,G_\alpha\left(\x,t+T(\x)\right)\qquad (v_\alpha >0)
\end{aligned}
\end{equation}
is a causal distribution i.e. $\mbox{supp}(g_\alpha(\x,\cdot))\subseteq [0,\infty)$. Here $T(\x) := \int_0^{|\x|} \frac{1}{v_\alpha(r)}\, \d r
\geq 0$ denotes the \emph{travel time} of $G_\alpha$. Then any attenuated wave $\A_\alpha(f)$ with compactly supported source $f$ has finite
front speed if and only if the Green function $G_\alpha(\x,\cdot)$ has finite front speed. We note that property~(\ref{defg}) implies
$$
   \mbox{supp}(G_\alpha(\x,\cdot)) \subseteq [0,\infty)\qquad\mbox{ for all $\x\in\R^3$}\,.
$$
Moreover, we assume that the front speed $v_\alpha$ is continuous, which implies $T(\x)\in C(\R^3)$ and $g_\alpha\in\mathcal{S}(\R,C(\R^3))$.
For the standard case were the wave front speed is assumed to be constant, this assumption is satisfied.

\subsection*{General structure of attenuated waves}

The property $g_\alpha\in\mathcal{S}(\R,C(\R^3))$ implies that
\begin{equation}\label{defg0}
\begin{aligned}
  g_\alpha(\x,t)
     &= \frac{1}{\sqrt{2\,\pi}}\,\F\left(e^{-\beta_*(|\x|,\omega)}\right)\,,\\
\end{aligned}
\end{equation}
with $\mbox{Re}(\beta_*):=-\log \left(\sqrt{2\,\pi}|\check g_\alpha|\right)$ and $\mbox{Im}(\beta_*):=-\mbox{arg}\left(\check
g_\alpha\right)$. Here $\check g_\alpha(\x,\cdot)$ denotes the inverse Fourier transform of $g_\alpha$. Since $g_\alpha$ is real-valued, it
follows that $\mbox{Re}(\beta_*)$ is even and $\mbox{Im}(\beta_*)$ is odd with respect to $\omega$. 
(We note that the causality of $g_\alpha(\x,\cdot)$ implies that $z\in\R+i\,\R^+\mapsto \check g_\alpha(\x,z)\in\C$ is analytic 
(cf. Theorem~\ref{th:hoer}) and thus $-\mbox{Re}(\beta_*(z)):\R+i\,\R^+\mapsto [-\infty,\infty)$ is subharmonic 
(cf. Example~4.1.10 in~\cite{Ho03}).) \\
If $\frac{\partial }{\partial r}\mbox{Im}\left(\beta_*(r,\omega)\right) <\frac{\omega}{v_\alpha(r)} $, then $\beta_*$ can be written 
as follows
\begin{equation}\label{defalpha0}
\begin{aligned}
  \beta_*(|\x|,\omega)
     &=
     a_S(\omega) +
     \int_0^{|\x|}\left[\alpha(r,\omega)
                 -i\, \left( \frac{\omega}{c(r,\omega)}-\frac{\omega}{v_\alpha(r)}  \right)\right]\, \d r \,,\\
\end{aligned}
\end{equation}
where $\alpha$ and $c$ are positive and even with respect to $\omega$ and $\lim_{\omega\to\infty} \mbox{Re}(\beta_*) = \infty$. 
Moreover, $a_S$ is a complex valued function such that
$$
%    g_S(t) := 
     g_\alpha(\mathbf{0},t)
    \in\mathcal{S}(\R)\,.
$$
By $S_\alpha$ we denote the time-convolution operator with kernel $g_\alpha(\mathbf{0},t)$ that maps $\mathcal{S}(\R)$ into 
$\mathcal{S}(\R)$. If $\beta_*$ depends only on 
$\omega$, i.e. $\beta_*=a_S$, then the attenuated wave satisfies the standard wave equation with source term $f$ replaced by 
$S_\alpha\,f$. We note that $g_\alpha$ is causal if and only if Theorem~\ref{th:hoer} (cf. Section~\ref{sec-causcheck}) is satisfied.

The standard normal form~(\ref{patt}) assumes that
$$
   S_\alpha=\mbox{Id},\,\quad v_\alpha=c_0=const,\,\quad \alpha=\alpha(\omega)\quad\mbox{ and }\quad 0<\alpha(\omega)<\infty\,,
$$
%the front speed is constant, say $c_0$, and that the attenuation law $\alpha$ depends only on $\omega$,
which together with the Kramers-Kronig relation~(\ref{KKst}) imply that $c=c(\omega)$. Moreover, 
$\beta_*(|\x|,\omega)=\alpha_*(\omega)\,|\x|$ with $\alpha_*$ defined as in~(\ref{defalpha*}). We note that $g_\alpha$ is causal if and only 
if $\alpha_*$ satisfies Lemma~\ref{lemm:caus1}.

Since wave attenuation is an irreversible thermodynamic process, $\beta_*$ can depend on time or
equivalently on the distance to the origin. We will see that in general, wave attenuation depends on its history, even if $\alpha$ depends
only on $\omega$. In particular this is the case if the support of the kernel of $S_\alpha$ is not discrete. Moreover, it is not evident that
the front speed of an attenuated spherical wave is constant. Therefore we see no reason to assume that $\beta_*$ and the
front speed $v_\alpha$ depend only on $\omega$. \\

%\noindent
\begin{rema}\label{rema:green}
If $v_\alpha(r)\leq v_B$  for all $r\geq 0$ ($v_B=const.$) and $g_\alpha$ is a causal distribution, then $\tilde
g(\x,t):=4\,\pi\,|\x|\,G_\alpha\left(\x,t+\frac{|\x|}{v_B}\right)$ must also be a causal distribution. In this case $\tilde \alpha_*$
corresponding to $\tilde g(\x,t)$ is given by~(\ref{defalpha0}) with $v_\alpha$ replaced by $v_B$ and the Green function reads as follows
$G_\alpha(\x,t)=\tilde g_\alpha(\x,t) *_t \frac{\delta(t-\frac{|\x|}{v_B})}{4\,\pi\,|\x|}$. This fact will be used to prove the non-causality
of some wave models.
\end{rema}

\section{Wave equation obeying attenuation and causality}
\label{sec-genwaveeq}

Now we derive the wave equation satisfied by the attenuated waves described in Section~\ref{sec-genmodel} and discuss its Cauchy problem.

%\subsection*{The derivation of the wave equation}

First we derive the wave equation for the Green function $G_\alpha$. The most convenient derivation uses the representation of the Green
function introduced in Remark~\ref{rema:green}:
\begin{equation}\label{pspalpha}
\begin{aligned}
    &\qquad\qquad\qquad G_\alpha = g_\alpha *_t G_B  \qquad\qquad\qquad\qquad\mbox{ with}\\
%\end{aligned}
%\end{equation}
%with
%\begin{equation}\label{Ginftyg}
%\begin{aligned}
    &G_B (\x,t) = \frac{\delta(t-\frac{|\x|}{v_B})}{4\,\pi\,|\x|}
    \quad\mbox{ and }\quad
    g_\alpha (\x,t) = \frac{1}{\sqrt{2\,\pi}}\,\F \left( e^{\beta_*(|\x|,\omega)} \right)(t)\,.
\end{aligned}
\end{equation}
Here the constant $v_B$ is an upper bound of the front speed of $G_\alpha$ and $\beta_*$ is defined as in~(\ref{defalpha0}) with
$v_\alpha(\x)$ replaced by $v_B$. We recall that $S_\alpha$ is the time convolution operator such that
$S_\alpha\,\delta(t)=g_\alpha(\mathbf{0},t)\in\mathcal{S}(\R)$. To formulate the wave equation we need the time convolution operators
$D_*:\mathcal{D}_+'(\R)\to \mathcal{D}_+'(\R)$ and $D_{*,|x|}:\mathcal{D}_+'(\R)\to \mathcal{D}_+'(\R)$ with causal kernels
\begin{equation}\label{defK*}
\begin{aligned}
   K_*(r,t) &= \frac{1}{\sqrt{2\,\pi}}\,\F \left( \frac{\partial\beta_*}{\partial r}(r,\omega) \right)(t)
   \quad\mbox{ and }\quad
   K_{*,|x|}(r,t) &= \frac{\partial K_*}{\partial r}(r,t)\,
\end{aligned}
\end{equation}
for all $r,\,t>0$, respectively. Here $\mathcal{D}_+'(\R)$ denotes the space of causal distributions. From~(\ref{pspalpha}), it follows
\begin{equation*}
\begin{aligned}
   \nabla g_\alpha &= -\frac{\x}{|\x|}\,K_* *_t g_\alpha \,\qquad\quad\mbox{ and }\\
   \nabla^2 g_\alpha &= \left[ -\frac{2}{|\x|}\,K_*
                              -K_{*,|x|}
                              + K_**_t K_*
                        \right] *_t g_\alpha\,.
\end{aligned}
\end{equation*}
This together with~(\ref{pspalpha}) imply 
\begin{equation*}
\begin{aligned}
   \nabla^2 G_\alpha
      -\frac{1}{v_B^2}\,\frac{\partial^2 G_\alpha}{\partial t^2}
       = \left[ D_*^2 +\frac{2}{v_B}\,\frac{\partial}{\partial t}\,D_* - D_{*,|x|}\right] G_\alpha
          - S_\alpha\,\delta(t)\,\delta(\x)\,.
\end{aligned}
\end{equation*}
Due to causality of $g_\alpha$ we have
\begin{equation}
\begin{aligned}
        \left.G_\alpha\right|_{t<0}=0\,,\quad
          \left.\frac{\partial G_\alpha}{\partial t}\right|_{t<0}=0\,.
\end{aligned}
\end{equation}
From
$$
   g_\alpha,\,\frac{\partial g_\alpha}{\partial x_j},\, \nabla^2 g_\alpha \in \mathcal{S}(\R,C(\R^3))\qquad \mbox{ for $j=1,\,2,\,3$}\,,
$$
it follows that
$$
 K_{\Lambda_\alpha}
   := \left[ D_*^2 +\frac{2}{v_B}\,\frac{\partial}{\partial t}\,D_* - D_{*,|x|}\right] G_\alpha \in \mathcal{S}(\R,C(\R^3))\,
$$
and that the space-time-convolution operator $\Lambda_\alpha$ with kernel $K_{\Lambda_\alpha}$ maps $\mathcal{S}(\R,C^\infty(\R^3))$ into
$\mathcal{S}(\R,C^\infty(\R^3))$. Since an arbitrary attenuated wave is of the form $p_\alpha=G_\alpha*_{\x,t}f$, the general wave equation
reads as follows:
\begin{equation}\label{genwaveeq}
\begin{aligned}
   &\nabla^2 p_\alpha
      -\frac{1}{v_B^2}\,\frac{\partial^2 p_\alpha}{\partial t^2}
       = (\Lambda_\alpha - S_\alpha)(f) \quad\mbox{ with }\\
        &\qquad\left.p_\alpha\right|_{t<0}=0\,,\quad
          \left.\frac{\partial p_\alpha}{\partial t}\right|_{t<0}=0\,.
\end{aligned}
\end{equation}
This equation has for every $f\in \mathcal{S}(\R,C^\infty(\R^3))$ with compact support a unique solution
\begin{equation*}
\begin{aligned}
   p_\alpha(\x,t) = \int_{\R^3} \frac{[\Lambda_\alpha(f)-S_\alpha(f)]\left(\x',t-\frac{|\x-\x'|}{v_B}\right)}{4\,\pi\,|\x-\x'|}  \,\d \x' \;\in\;
   \mathcal{S}(\R,C^\infty(\R^3))
\end{aligned}
\end{equation*}
with finite wave front speed.  \\

If the attenuation law and the phase speed do not depend on the spatial position, then the operator $ D_*$ does not depend on the spatial
position, too, and
\begin{equation}\label{D*rel}
D_*(G_\alpha)*_{\x,t} f =  D_*(p_\alpha)
\end{equation}
holds. Moreover, then $D_{*,|x|}$ is the zero operator. In this case we can write the wave equation as follows
\begin{equation}\label{genwaveeq0}
\begin{aligned}
   &\nabla^2 p_\alpha
      -\left[ D_* +\frac{1}{v_B}\,\frac{\partial }{\partial t} \right]^2 p_\alpha
       =  - S_\alpha(f)\,,
%\\
%       & \left.p_\alpha\right|_{t<0}=0\,,\quad
%          \left.\frac{\partial p_\alpha}{\partial t}\right|_{t<0}=0\,.
\end{aligned}
\end{equation}
with $\left.p_\alpha\right|_{t<0}=0$ and $\left.\frac{\partial p_\alpha}{\partial t}\right|_{t<0}=0$. In general the supports of the kernels
of $D_*$ and $S_\alpha$ are subsets of $[0,\infty)$ with positive Lebesgue measure, which means that the attenuated wave depends on its
history. Since the values of the wave in the past are required, the Cauchy problem of wave equation~(\ref{genwaveeq0}) is not reasonable. In
the next theorem we formulate a generalization of the Cauchy problem and state its properties for the special case $\alpha=\alpha(\omega)$,
$c=c(\omega)$ with $a_S=0$.

\begin{prop}
Let $D_*$ be the time-convolution operator with causal kernel $K_*$ defined as in~(\ref{defK*}) and let
$q\in\mathcal{S}(\R^3\times\R_0^-)$, $\varphi:=\left.q\right|_{t=0}$ and $\psi:=\left.\frac{\partial q}{\partial t}\right|_{t=0}$. Moreover,
let $G_\alpha$ denote the Green function of wave equation~(\ref{genwaveeq0}) with $S_\alpha=\mbox{Id}$ and wave front speed $\leq v_B$.
Then the solution of the generalized Cauchy problem
\begin{equation*}
\begin{aligned}
   &\nabla^2 p_\alpha
      -\left[ D_* + \frac{1}{v_B}\,\frac{\partial }{\partial t}\right]^2\,p_\alpha = 0
      \qquad\mbox{ on $\R^3\times \R_0^+$}\,,\\
       & \left.p_\alpha\right|_{t\leq 0}=q
       \quad\mbox{ and }\quad
          \left.\frac{\partial p_\alpha}{\partial t}\right|_{t=0}=\psi\,,
\end{aligned}
\end{equation*}
is given by
\begin{equation*}
\begin{aligned}
    p_\alpha = G_\alpha *_{\x,t} f
    \qquad\mbox{ on $\R^3\times \R_0^+$}\,,
\end{aligned}
\end{equation*}
with
\begin{equation*}
\begin{aligned}
      f = -  \frac{\psi}{v_B^2}\,\delta(t)
              -  \frac{\varphi}{v_B^2} \,\delta'(t)
                 - \left[ D_*^2 + \frac{1}{v_B}\,\frac{\partial }{\partial t}\,D_* , M_{H(t)}\right]\,q \,.
\end{aligned}
\end{equation*}
Here $M_{H(t)}$ denotes the multiplication operator, $H=H(t)$ the Heaviside function and $[,]$ denotes the commutator, i.e. if $A$, $B$ are
operators then $[A,B]=A\,B-B\,A$.
\end{prop}

\begin{proof}
Let $p_\alpha$, $G_\alpha$, $f$ and $H(t)$ be defined as in the Proposition. For convenience  let $A_*:=D_*^2 +
\frac{1}{v_B}\,\frac{\partial}{\partial t}\,D_*$. Then $\tilde p_\alpha:= H(t)\,p_\alpha$ satisfies the following properties:
$$
\nabla^2 \tilde p_\alpha = H(t)\,\nabla^2 p_\alpha\,,
$$
\begin{equation*}
\begin{aligned}
   \frac{\partial^2 \tilde p_\alpha}{\partial t^2}
   &=  H(t) \,\frac{\partial^2 p_\alpha}{\partial t^2}
        +  \psi\,\delta(t)
              + \varphi \,\delta'(t)\,,
\end{aligned}
\end{equation*}
and
\begin{equation*}
\begin{aligned}
  A_* \,\tilde p_\alpha
   = H\,A_*(p_\alpha) + [ A_*,M_H]\,p_\alpha\,,
\end{aligned}
\end{equation*}
since $\varphi=\left.p_\alpha\right|_{t=0}$ and $\psi=\left.\frac{\partial p_\alpha}{\partial t}\right|_{t=0}$. From these properties we infer
\begin{equation*}
\begin{aligned}
   &\nabla^2 \tilde p_\alpha
      -\frac{1}{v_B^2}\,\frac{\partial^2 \tilde p_\alpha }{\partial t^2} - A_*\,\tilde p_\alpha
      =
      - \frac{\psi}{v_B^2}\,\delta(t)
      - \frac{\varphi}{v_B^2} \,\delta'(t)
      - [ A_*,M_H]\,p_\alpha\,,\\
       & \left.\tilde p_\alpha\right|_{t< 0}=0
       \quad\mbox{ and }\quad
          \left.\frac{\partial \tilde p_\alpha}{\partial t}\right|_{t<0}=0\,,
\end{aligned}
\end{equation*}
which has the solution $\tilde p_\alpha=G_\alpha*_{\x,t} f$ and thus $p_\alpha=G_\alpha*_{\x,t} f$ on $\R^3\times \R_0^+$. This proves the
Proposition.
\end{proof}

\section{Causality analysis of attenuated spherical waves obeying the frequency power law}
\label{sec-causcheck}

In this section we show that not every attenuated spherical wave of the form~(\ref{patt}) has a finite front speed although its
phase speed is properly related to the frequency power law via the Kramers-Kronig relation~(\ref{KKst}).
Other attenuated spherical wave models are analysed in Section~\ref{sec-szabo},~\ref{sec-caustv} and ~\ref{sec-gamgen} \\

\noindent In the following we use the notions $\R^+=(0,\infty)$, $\R_0^+=[0,\infty)$, $\R^-=(-\infty,0)$ and
$\R_0^-=(-\infty,0]$.
The next Theorem is a reformulation of Theorem~7.4.3 in~\cite{Ho03} for the case of causal tempered distributions.
%We recall that a tempered distribution $u\in \mathcal{S}'(\R)$ is called causal, if its support is contained in $\R_0^+$.
\begin{theo}\label{th:hoer}
A distribution $u\in\mathcal{S}'(\R)$ is causal, i.e. $\mbox{supp}(u)\subseteq \R_0^+$,
if and only if
\begin{description}
\item [(A1)] $\R+i\,\R^-\to\C,\,z \mapsto \hat u(z)$ is analytic and\\
\item [(A2)] $\exists \epsilon>0\,\exists C>0\,\exists N>0\,\forall z\in \R+i\,(-\infty,-\epsilon):\,|\hat u(z)|\leq C\,(1+|z|)^N\,.$
\end{description}
%If items 1) and 2) are satisfied then
%$\F(u(t)\,e^{\eta\,t})(\omega)=\hat u(\omega+i\,\eta)$ for all $\eta\in\R^-$.
\end{theo}

Applying Theorem~\ref{th:hoer} to attenuated spherical waves of the form~(\ref{patt}) yields the
following Lemma.

\begin{lemm}\label{lemm:caus1}
Let $p_\alpha$ be defined by~(\ref{patt}) with real-valued functions $\alpha=\alpha(\omega)>0$ and $c=c(\omega)$,
and let $p_\alpha(\x,\cdot)\in\mathcal{S}'(\R)$ for any $\x\in\R^3$.
%Let $\alpha=\alpha(\omega)>0$ and $c=c(\omega)$ be real-valued functions on $\R$.
The wave defined by~(\ref{patt}) has a finite
wave front speed if and only if $\alpha_*$ defined by~(\ref{defalpha*}) satisfies the following conditions:
\begin{description}
\item [(B1)] $\alpha_*=\alpha_*(-z)$ is analytic on $\R+i\,\R^-$ and\\
\item [(B2)] $\exists \epsilon>0\,\exists C>0\,\exists N>0\,\forall z\in \R+i\,(-\infty,-\epsilon):$\\
             $\mbox{\;\;\;\;\;\;\;\;\;\;\;\;\;\;} -\mbox{Re}(\alpha_*(-z)) \leq C + N\,\log (1+|z|)\,.$
\end{description}
\end{lemm}

\begin{proof}
According to Theorem~\ref{th:hoer} $g_\alpha(\x,\cdot)$ defined as in~(\ref{defg}) is a causal distribution for any fixed $\x\in\R^3$, if
properties (A1) and (A2) hold for $u:=g_\alpha(\x,\cdot)$ for any $\x\in\R^3$ . Let $\x\in\R^3$ be arbitrary but fixed.
According to~(\ref{Fst}) we have $g_\alpha(\x,\cdot) =
\F\left(e^{-\alpha_*(\omega)\,|\x|}\right)=\F^{-1}\left(e^{-\alpha_*(-\omega)\,|\x|}\right)$ and thus
\begin{equation}\label{Fhat}
\begin{aligned}
    \hat g_\alpha(\x,\omega) = e^{-\alpha_*(-\omega)\,|\x|} \qquad(\omega\in\R)\,.
\end{aligned}
\end{equation}
Let $\beta_1(z):=-\mbox{Re}(\alpha_*(-z))\,|\x|$ and $\beta_2(z):=-\mbox{Im}(\alpha_*(-z))\,|\x|$.
Since $e^{-\alpha_*(-z)\,|\x|}$ is analytic on $\R+i\,\R^-$, the Cauchy-Riemann equations are satisfied
\begin{equation}\label{creq}
\begin{aligned}
    \left[\frac{\partial \beta_1}{\partial x} - \frac{\partial \beta_2}{\partial y}\right]\,\cos\beta_2
      &= \left[\frac{\partial \beta_1}{\partial y} + \frac{\partial \beta_2}{\partial x}\right]\,\sin\beta_2\,,\\
 \left[\frac{\partial \beta_1}{\partial x} - \frac{\partial \beta_2}{\partial y}\right]\,\sin\beta_2
     &= -\left[\frac{\partial \beta_1}{\partial y} + \frac{\partial \beta_2}{\partial x}\right]\,\cos\beta_2\,,
\end{aligned}
\end{equation}
which imply for all $z$ with $\cos\beta_2(z)\not=0$ and $\sin\beta_2(z)\not=0$  the equations
\begin{equation*}
\begin{aligned}
    \frac{\partial \beta_1}{\partial x} = \frac{\partial \beta_2}{\partial y}
\qquad\mbox{ and }\qquad
       \frac{\partial \beta_1}{\partial y} = - \frac{\partial \beta_2}{\partial x}\,.
\end{aligned}
\end{equation*}
This means that $\beta_1(z)+i\,\beta_2(z)$ is analytic for all  $z\in\R+i\,\R^-$ satisfying $\cos\beta_2(z)\not=0$ and $\sin\beta_2(z)\not=0$.
The same equations follow easily if $\sin\beta_2(z)\not=0$ but $\cos\beta_2(z)=0$, and for $\cos\beta_2(z)\not=0$ but $\sin\beta_2(z)=0$ and
thus $\alpha_*(-z)$ is analytic on  $\R+i\,\R^-$. This shows that condition (B1) is satisfied. Conversely if $\alpha_*(-z)$ is analytic on
$\R+i\,\R^-$ then due to the chain rule $e^{-\alpha_*(-z)\,|\x|}$ must be analytic on $\R+i\,\R^-$, since the complex exponential function is
analytic on $\C$.

From condition (A2) in Theorem~\ref{th:hoer} together with~(\ref{Fhat}) we infer
$$
%    \exists C>0\,\exists N>0\,\forall z\in \R-i\,(0,\infty):\,
      e^{-\mbox{Re}(\alpha_*(-z)\,|\x|)} \leq C\,(1+|z|)^N\,,
$$
whereby we can assume without loss of generality $C>1$. Let $\tilde C:=\log C>0$. Since  the real (natural) logarithm function is monotonic
increasing, we can apply it onto the previous inequality, which yields condition (B2) with $C$ replaced by $\frac{\tilde C}{|\x|}$.
Conversely, if (B2) holds then condition (A2) holds, too. This concludes the proof.
\end{proof}

For the analysis of the causality properties of the standard models of spherical attenuated waves we need the
following Lemma.

\begin{lemm} \label{lemm:caus2}
Let $\x\in\R^3$ and $\gamma>0$ be arbitrary but fixed. Moreover, let $s(\gamma)$ be the sign of $\cos(\frac{\pi}{2}\,\gamma)$
if $\cos(\frac{\pi}{2}\,\gamma)\not=0$ and let $s(\gamma)=1$ if $\cos(\frac{\pi}{2}\,\gamma)=0$. 
The distribution $\F((-i\,\omega)^\gamma)$ is causal on $\R$ for any $\gamma>0$ and
the distribution $\F(e^{-s(\gamma)\,(-i\,\omega)^\gamma\,|\x|})$ is causal if $\gamma\in (0,1]$
and non-causal if $\gamma>1$.
\end{lemm}

\begin{proof}
i) First we prove the causality of $\F((-i\,\omega)^\gamma)$.
Since $\omega\in \R\mapsto (-i\,\omega)^\gamma\in\C$ is a slowly increasing function, it is a tempered
distribution and hence $t\in\R\mapsto\F((-i\,\omega)^\gamma)(t)\in\R$ is also a tempered distribution.
Therefore Theorem~\ref{th:hoer} can be applied.
Let $\gamma\in (0,\infty)$ be arbitrary but fixed and let $U$ denote the complex plane without the
non-negative real axis. Then $z^\gamma:=e^{\gamma\,\log z}$ is analytic on U (cf.\cite{La93}).
Therefore $(i\,z)^\gamma$ is analytic on $\R+i\,\R^-$ and condition (A1) in Theorem~\ref{th:hoer} is satisfied.
Moreover, $z\in U\mapsto (i\,z)^\gamma\in\C$ satisfies condition (A2) in Theorem~\ref{th:hoer} with $C=1$ and
$N=\gamma$, which shows that $\F((-i\,\omega)^\gamma)(t)$ is causal on $\R$.

ii) For the rest of the proof let $\x\in\R^3$ be arbitrary but fixed.
Now we prove that $\sqrt{2\,\pi}\,g_\alpha(\x,t):=\F(e^{-s(\gamma)\,(-i\,\omega)^\gamma\,|\x|})$ is causal for $\gamma\in (0,1]$ and
non-causal for $\gamma\in (1,\infty)$.
Since $(-i\,\omega)^\gamma = |\omega|^\gamma\,[\cos(\frac{\pi}{2}\,\gamma)-i\,\mbox{sgn}(\omega)\,\sin(\frac{\pi}{2}\,\gamma)]$, it follows that 
$\left|e^{-s(\gamma)\,(-i\,\omega)^\gamma\,|\x|}\right|=e^{-|\cos(\frac{\pi}{2}\,\gamma)|\,|\omega|^\gamma\,|\x|}$ is bounded
for each $\gamma\in \R^+$ and thus $\F(e^{-s(\gamma)\,(-i\,\omega)^\gamma\,|\x|})$ is a tempered distribution
and Lemma~\ref{lemm:caus1} can be applied. Using the same notion as in~(\ref{defg0}), yields 
$\alpha_*(-z)=s(\gamma)\,(i\,z)^\gamma$. \\
Let $z=r\,e^{i\,\varphi}$ with $r>0$ and $\varphi\in (-\pi,0)$. Then $z\in\R+i\,\R^-$ and
$\mbox{Re}((i\,z)^\gamma) = \cos(\gamma(\varphi+\frac{\pi}{2}))\,|z|^\gamma$, and thus the inequality in
(B2) reads as follows
\begin{equation}\label{lem:help}
     -s(\gamma)\,\cos(\gamma(\varphi+\frac{\pi}{2}))\,\,|z|^\gamma   \leq   C + N\,\log (1+|z|)
\qquad\quad (\gamma>0)\,.
\end{equation}
This shows that condition (B2) is satisfied if and only if $s(\gamma)\,\cos(\gamma(\varphi+\frac{\pi}{2}))\geq 0$. \\
a) For $\gamma\in (0,1]$ we get $s(\gamma)=1$ and $\gamma(\varphi+\frac{\pi}{2})\in (-\frac{\pi}{2},\frac{\pi}{2})$
for any $\varphi\in (-\pi,0)$ and thus $s(\gamma)\,\cos(\gamma(\varphi+\frac{\pi}{2}))\geq 0$, i.e.
condition (B2) is satisfied for any $\epsilon>0$. Therefore $\F(e^{-s(\gamma)\,(-i\,\omega)^\gamma\,|\x|})$ is causal for $\gamma\in (0,1]$.\\
b) Now we prove the non-causality of $\F(e^{-s(\gamma)\,(-i\,\omega)^\gamma\,|\x|})$ for
$\gamma\in [3,5]\cup [7,9]\cup\cdots$. Since for these $\gamma-$values $s(\gamma)=1$ we have
to find a $\varphi$ such that the sign of $\cos(\gamma(\varphi+\frac{\pi}{2}))$ is negative.
For $\gamma>1$ let $0 < \delta<\mbox{min}\left(\frac{\pi}{4},\frac{\gamma-1}{\gamma+1}\,\frac{\pi}{2}\right)$
and $\varphi_\delta:=\left(\frac{1}{\gamma}-1\right)\,\frac{\pi}{2}+\frac{\delta}{\gamma}$.
Since $\varphi_\delta\in (-\frac{\pi}{2},-\delta)$ and
$\gamma(\varphi_\delta+\frac{\pi}{2})=\frac{\pi}{2}+\delta\in (\frac{\pi}{2},\frac{3\pi}{4})$,
it follows that $z:=r\,e^{i\,\varphi_\delta}\in\R+i\,\R^-$ and
$\cos(\gamma(\varphi_\delta+\frac{\pi}{2}))<0$ and hence condition (B2) cannot be satisfied for any $\epsilon>0$. \\
c)  Now we prove the non-causality of $\F(e^{-s(\gamma)\,(-i\,\omega)^\gamma\,|\x|})$ for
$\gamma\in (1,3)\cup (5,7)\cup (9,11)\cup\cdots$. Since for these $\gamma-$values $s(\gamma)=-1$
we have to find a $\varphi$ such that the sign of $\cos(\gamma(\varphi+\frac{\pi}{2}))$ is positive.
For $\varphi=-\frac{\pi}{2}$ it follows that $z=r\,e^{i\,\varphi}\in i\,\R^-$ ($r>0$) and
$\cos(\gamma(\varphi+\frac{\pi}{2}))=1$, which implies at once that condition (B2) cannot be satisfied for any $\epsilon>0$.
This proves the Corollary.
\end{proof}

The following two Corollaries clarify for which values of $\gamma$ and $\omega_0$ the frequency power law $\alpha(\omega) =
\alpha_0\,|\omega|^\gamma$ together with the phase speed determined by the Kramers-Kronig relations~(\ref{KKst}), yield an attenuated
spherical wave that satisfies the causality requirement. The derivation of the phase speed $c=c(\omega)$ corresponding to frequency power law
with various $\gamma-$values can be found in~\cite{Wa00,WaMoMi00,WaHuMoMi03,WaHuMoBrMi00,Szabo95}.

\begin{coro} \label{coro:powlaw1}
Let $\alpha_0$ and $\omega_0$ be positive constants, $\gamma\in \R^+\backslash\N$ and  let
the attenuation law $\alpha$ and the phase speed $c$ be defined as follows:
\begin{equation}\label{alpha1}
\begin{aligned}
    \alpha(\omega) &= \alpha_0\,|\omega|^\gamma\,,\\
   \frac{1}{c(\omega)} - \frac{1}{c(\omega_0)} &= \alpha_0\,\tan\left(\frac{\pi}{2}\,\gamma\right)\,
                                   (|\omega|^{\gamma-1} - |\omega_0|^{\gamma-1})\,.
\end{aligned}
\end{equation}
Then $\alpha_*$ defined as in~(\ref{defalpha*}) reads as follows
\begin{equation}\label{beta1}
\begin{aligned}
\alpha_*(\omega)
    = \frac{\alpha_0\,(-i\omega)^\gamma}{\cos(\frac{\pi}{2}\,\gamma)}
       + i\alpha_0\, \tan\left(\frac{\pi}{2}\,\gamma\right)\,|\omega_0|^{\gamma-1} \,\omega
\end{aligned}
\end{equation}
and the wave defined by~(\ref{patt}) has finite front speed if $\gamma\in (0,1)$ and $\omega_0=0$. For these cases the wave front speed is
equal to $c(0)$. If $\omega_0=0$ and $\gamma\in (1,\infty)\backslash\N$ or $\omega_0>0$ and $\gamma\in\R^+\backslash\N$, then the wave 
defined by~(\ref{patt}) cannot have a finite front speed.
\end{coro}

\begin{proof}
Let $\x\in\R^3$ be arbitrary but fixed.  From~(\ref{alpha1}) and~(\ref{defalpha*}) we get at once identity~(\ref{beta1}). \\
i) Let $\omega_0=0$, then according to Lemma~\ref{lemm:caus2} $\F(e^{-\alpha_*\,|\x|})(t)$ with $\alpha_*$ defined as in~(\ref{beta1})
is causal if $\gamma\in (0,1)$ and non-causal if $\gamma\in (1,\infty)\backslash\N$. \\
ii) Now let $\omega_0>0$. Then it follows that
\begin{equation*}
  -\mbox{Re}\left( \alpha_*(-z)\right)
      = - a_1(\gamma)\,|z|^\gamma + a_2(\gamma)\,|\mbox{Im}(z)|
\qquad\mbox{ for all $z\in\R+i\,\R^-$}\,,
\end{equation*}
where $a_1(\gamma):=\alpha_0\,\frac{\cos((\varphi+\frac{\pi}{2})\,\gamma)}{\cos(\frac{\pi}{2}\,\gamma)}$, $\varphi\in (-\pi,0)$ is the
argument of $z$ and $a_2(\gamma):=\alpha_0\, \tan\left(\frac{\pi}{2}\,\gamma\right)\,|\omega_0|^{\gamma-1}$.\\
a) If $\gamma\in (0,1)$ then $a_1,\,a_2>0$ and
\begin{equation*}
   -\mbox{Re}\left( \alpha_*(-i\,z_2)\right) =
       - a_1(\gamma)\,|z_2|^\gamma + a_2(\gamma)\,|z_2|
\qquad\mbox{ for all $z_2<0$}\,,
\end{equation*}
i.e. $-\mbox{Re}\left( \alpha_*(-i\,z_2)\right)$ growths like $|z_2|$ for sufficiently large $-z_2$. Since this term cannot be bounded by
$\log(1+|z_2|)$, condition (B2) in Lemma~\ref{lemm:caus1} cannot be satisfied for $\omega_0>0$ and $\gamma\in (0,1)$. \\
b) Now let $\gamma>1$ and $\gamma\not\in\N$. We note that $a_1(\gamma)$ has the same sign as 
$s(\gamma)\,\cos((\varphi+\frac{\pi}{2})\,\gamma)$, where $s(\gamma)$ is defined as in Lemma~\ref{lemm:caus2}. 
As in the proof of Lemma~\ref{lemm:caus2} one shows that for an appropriate choice of $\varphi\in (-\pi,0)$
the constant $a_1$  is negative, which shows that $-\mbox{Re}\left(\alpha_*(-z(r))\right)$ growth like $|z|^\gamma$ for sufficiently
large $|z|$. Therefore condition (B2) in Lemma~\ref{lemm:caus1} cannot be
satisfied for $\omega_0>0$ and $\gamma\in (0,1)$. In summary we have shown that the front speed of the wave defined by~(\ref{patt})
and~(\ref{beta1}) cannot be finite if $\omega_0>0$ and $\gamma\in\R^+\backslash\N$. \\
iii) Now we show that the front speed for the case $\omega_0=0$ with $\gamma\in (0,1)$ is equal to $c(0)$.
We recall that if $g_\alpha(\x,t)$ is causal, then the front wave speed satisfies
the condition $v_\alpha\leq c(\omega_0)$. If the front wave speed $v_\alpha$ at $\x$ is smaller than
$c(0)$ ($\omega_0=0$), then $g_\alpha(\x,t+\delta\,|\x|)$ is causal for some $\delta> 0$.
This means that there exist constants $\epsilon>0$, $C>0$ and $N>0$ such that for all $z\in\R+i\,(-\infty,-\epsilon)$:
\begin{equation}\label{coro:help1}
    -\mbox{Re}(\alpha_*(-z)) - \mbox{Re}(i\,(-z)\,\delta) \leq C + N\,\log (1+|z|)\,.
\end{equation}
Since $\omega_0=0$, we have $a_2=0$. As above we get for all $z=i\,z_2\in\R^-$
\begin{equation*}
  -\mbox{Re}\left( \alpha_*(-z)\right)
      = - a_1\,|z_2|^\gamma + \delta\,|z_2|  \qquad (a_1>0)\,,
\end{equation*}
which grows like $|z_2|$ if $\delta>0$ and thus condition~(\ref{coro:help1}) cannot be satisfied. This contradicts the fact that
$g_\alpha(\x,t)$ is causal and hence we conclude that $\delta=0$. This proves that the wave front speed is equal to $c(0)$ and concludes the
proof.
\end{proof}

Now we come to the special case $\gamma=1$.

\begin{coro} \label{coro:powlaw2}
Let $\alpha_0,\,\omega_0>0$ be constants and
\begin{equation}\label{defbeta2}
\begin{aligned}
   \alpha_*(\omega)
      := \lim_{\gamma\to 1} \left[
              \frac{\alpha_0\,(-i\omega)^\gamma}{\cos(\frac{\pi}{2}\,\gamma)}+i\alpha_0\,
                 \tan\left(\frac{\pi}{2}\,\gamma\right)\,|\omega_0|^{\gamma-1} \,\omega
               \right]\,,
\end{aligned}
\end{equation}
then
\begin{equation}\label{betaalpha2}
\begin{aligned}
    \alpha(\omega) = \alpha_0\,|\omega|\,,
\qquad\mbox{ }\qquad
   \frac{1}{c(\omega)} - \frac{1}{c(\omega_0)} = -\alpha_0\,\frac{2}{\pi}\,
                                \log \left|\frac{\omega}{\omega_0}\right|\,
\end{aligned}
\end{equation}
and the wave defined by~(\ref{patt}) cannot have a finite front speed. Moroever, $\F(\alpha_*(\omega))$ is not causal.
\end{coro}

\begin{proof}
Definition (\ref{defbeta2}) implies
\begin{equation*}
\begin{aligned}
   \alpha_*(\omega)
      = \alpha_0\,|\omega|
            - i\alpha_0 \,\omega\,\lim_{\epsilon\to 0+}
                 \frac{|\omega|^\epsilon - |\omega_0|^\epsilon}{\cot(\frac{\pi}{2}\,(1+\epsilon))}\,.
\end{aligned}
\end{equation*}
Since both the numerator and the denominator of the last expression converge to zero, we can apply the \emph{rule of de l'Hospital}, 
which yields
\begin{equation*}
\begin{aligned}
   \alpha_*(\omega)
      -\alpha_0\,|\omega|
            &= -i\alpha_0 \,\omega\,\lim_{\epsilon\to 0+}
                 \frac{|\omega|^\epsilon\,\log|\omega|  - |\omega_0|^\epsilon\,\log|\omega_0|}
                      {-\frac{\pi}{2}\,\sin^2(\frac{\pi}{2}\,(1+\epsilon))} \\
            &= i\alpha_0\,\frac{2}{\pi}\,\omega\,(\log|\omega|-\log|\omega_0|)\,.
\end{aligned}
\end{equation*}
The result for the limit $\epsilon\to 0-$ follows analogously.
Comparing $\alpha_*$ with~(\ref{defalpha*}) yields~(\ref{betaalpha2}).
For $(x,y)\in\R^2$ let
$$
 u(x,y):= \alpha_0\,\sqrt{x^2+y^2} -\alpha_0\,\frac{2}{\pi}\,y\,(\log \sqrt{x^2+y^2}-\log|\omega_0|)
$$
and
$$
 v(x,y):= \alpha_0\,\frac{2}{\pi}\,x\,(\log \sqrt{x^2+y^2}-\log|\omega_0|)\,.
$$
Then $\alpha_*(z)=u(x,y)+i\,v(x,y)$ for $z\in\C$ with $x=\mbox{Re}(z)$ and $y=\mbox{Im}(z)$, and
$\frac{\partial u(x,y)}{\partial x}\not=\frac{\partial v(x,y)}{\partial y}$ for every $(x,y)\in\R^2\backslash {(0,0)}$.
Since $u(x,y)$ and $v(x,y)$ do not satisfy the Cauchy-Riemann equations for every $(x,y)\in\R^2\backslash {(0,0)}$,
$\alpha_*(-z)$ is not analytic on $\R+i\,\R^-$, i.e. condition (B1) in Lemma~\ref{lemm:caus1} is not satisfied.
This shows that for fixed $\x\in\R^3$ $\F(\alpha_*(\omega))(t)$ and
$\F(e^{-\alpha_*(\omega)\,|\x|})(t)$ cannot be causal and concludes the proof.
\end{proof}

\section{Causality analysis of Szabo's wave equation}
\label{sec-szabo}

The Green function of Szabo's equation~(\ref{szaboseq}) for $\gamma>0,\,\gamma\not\in\N$ is given by~(\ref{patt})  with the
following attenuation law and phase speed:
\begin{equation*}
\begin{aligned}
     \alpha(\omega) &= \mbox{Re}(\alpha_*(\omega))
\qquad\mbox{ and }\qquad
    \frac{1}{c(\omega)} - \frac{1}{c_0} = -\frac{\mbox{Im}(\alpha_*(\omega))}{\omega}  \,,
\end{aligned}
\end{equation*}
where 
\begin{equation}\label{alpha*szabo}
\begin{aligned}
    \alpha_*(\omega) &= \frac{1}{c_0}\,\sqrt{(-i\,\omega)^2
                          + 2\,\alpha_0\,c_0\,\frac{(-i\,\omega)^{\gamma+1}}{\cos\left(\frac{\pi}{2}\,\gamma\right)}}
                          + i\,\frac{\omega}{c_0}\,
\end{aligned}
\end{equation}
(cf. dispersion relation above equation~(\ref{szaboseq}).)
Here the square root is understood as the primitive square root, since $\alpha(\omega)$ has to be positive.
%If $\gamma\in (0,1)$ then the square root is nothing else but the primitive square root.
Since for $\gamma\in (0,1)$
%then for sufficiently large frequencies we have
\begin{equation*}
\begin{aligned}
  \alpha_*(\omega)
         &\approx \frac{(-i\,\omega)^{\gamma}}{\cos\left(\frac{\pi}{2}\,\gamma\right)}
      \qquad\mbox{for $|\omega|>>1$}
\end{aligned}
\end{equation*}
and for $\gamma>1$, $\gamma\not\in\N$
%for sufficiently small frequencies we have
\begin{equation*}
\begin{aligned}
  \alpha_*(\omega)
         &\approx \frac{(-i\,\omega)^{\gamma}}{\cos\left(\frac{\pi}{2}\,\gamma\right)}
      \qquad\mbox{ for $|\omega|<<1$}\,,
\end{aligned}
\end{equation*}
Szabo's model is a high frequency approximation of equation~(\ref{standwaveeq}) if $\gamma \in (0,1)$
and as a small frequency approximation of equation~(\ref{standwaveeq}) if $\gamma >1$, $\gamma\not\in\N$.

The next Proposition investigates the causality behaviour of Szabo's equation.

\begin{prop}\label{theo:szabo}
Let $\gamma\in \R^+\backslash\N$.
The Green function $G_\alpha$ of Szabo's equation~(\ref{szaboseq}) has finite front speed only if $\gamma\in (0,1)$. For $\gamma\in (0,1)$,  
the front speed of $G_\alpha$ is $c_0$.
\end{prop}

\begin{proof}
For the proof let $\x\in\R^3$ be arbitrary but fixed. The Green function $G(\x,t)$ of wave equation~(\ref{szaboseq}) has a front speed $\leq
c_0<\infty$, if $G(\x,t+\frac{|\x|}{c_0})$ is causal, i.e. if $\alpha_*(\omega)$ defined by~(\ref{alpha*szabo}) satisfies
Lemma~\ref{lemm:caus1}.
For convenience we set $\tilde \alpha_0:=\frac{\alpha_0}{\cos\left(\frac{\pi}{2}\,\gamma\right)}$. \\
i) First we prove the Proposition for $\gamma\in (0,1)$. 
Since $\sqrt{(i\,z)^2 + 2\,\tilde\alpha_0\,c_0\,(i\,z)^{\gamma+1}}$ (primitive square root) maps $\R+i\,\R^-$ analytically into 
$\R+i\,\R^-$, $\alpha_*(-z)$ maps $\R+i\,\R^-$ analytically into $\R+i\,\R^-$. This proves condition (B1). \\
Let $z=z_1+i\,z_2$ with $z_1\in\R$, $z_2\in\R^-$ and  
$$
   B(z) := \sqrt{1 + 2\,\tilde\alpha_0\,c_0\, \left(\frac{1}{i\,z}\right)^{1-\gamma}} - 1\,,
$$
then $\mbox{Re}\left(\alpha_*(-z)\right)=\mbox{Re}\left(i\,z\,B(z)\right)= -z_1\,\mbox{Im}\left(B(z)\right)-z_2\,\mbox{Re}\left(B(z)\right)$. 
Moreover, for Proposition~\ref{theo:szabo} property (B2) in Lemma~\ref{lemm:caus1} reads as follows: there exist constants
$\epsilon>0$, $C>0$ and $N>0$ such that for all $z\in \R+i\,(-\infty,-\epsilon)$:
\begin{equation}\label{rootszabo}
\begin{aligned} 
    \frac{1}{c_0}\,\left( z_1\,\mbox{Im}\left(B(z)\right) + z_2\,\mbox{Re}\left(B(z)\right) \right) 
%    -\frac{1}{c_0}\,\mbox{Re}\left(\sqrt{(i\,z)^2 + 2\,\tilde\alpha_0\,c_0\,(i\,z)^{\gamma+1}} - i\,z   \right)
    \leq C + N\,\log (1+|z|)\,.
\end{aligned}
\end{equation}
We prove this inequality by showing that its left hand side is always negative or zero.  \\
a) First we  prove that $\mbox{Re}\left(B(z)\right)\geq 0$ if $z_2\in\R^-$, which implies 
$z_2\,\mbox{Re}\left(B(z)\right)\leq 0$ if $z_2\in\R^-$. 
From $\gamma\in (0,1)$ and $z_2<0$, it follows that $(i\,\,z)^{1-\gamma}$ has positive real part. Since the inversion of a complex 
number with positive real part yields a complex number with positive real part, 
$1 + 2\,\tilde\alpha_0\,c_0\, \left(\frac{1}{i\,z}\right)^{1-\gamma}$ has a real part $\geq 1$. This implies 
that $B(z)$ has positive real part.\\
b) Now we show that $z_1\,\mbox{Im}\left(B(z)\right)\leq 0$ if $z_2\in\R^-$ and $z_1\in\R$. Let $z_1>0$, then the imaginary part of 
$(i\,z)^{1-\gamma}$ is positive and since the inversion of a complex number with positive imaginary part yields a complex number with 
negative imaginary part, $1 + 2\,\tilde\alpha_0\,c_0\, \left(\frac{1}{i\,z}\right)^{1-\gamma}$ has negative imaginary part. This implies 
that $B(z)$ has negative imaginary part and therefore $z_1\,\mbox{Im}\left(B(z)\right)\leq 0$ if $z_2\in\R^-$ and $z_1>0$. 
Now let $z_1<0$. Then the imaginary part of $(i\,z)^{1-\gamma}$ is negative and since the inversion of a complex number with negative 
imaginary part yields a complex number with positive imaginary part, $1 + 2\,\tilde\alpha_0\,c_0\, \left(\frac{1}{i\,z}\right)^{1-\gamma}$ 
has positive imaginary part. This implies that $B(z)$ has positive imaginary part and therefore $z_1\,\mbox{Im}\left(B(z)\right)\leq 0$ 
if $z_2\in\R^-$ and $z_1<0$. Clearly, if $z_1=0$ then $z_1\,\mbox{Im}\left(B(z)\right)= 0$. 
In summary we have proven that the left hand side of~(\ref{rootszabo}) is smaller or equal to zero and thus the inequality holds.\\
ii) The second part of the theorem is first proven for $\gamma\in (1,3)\backslash\{2\}$ and then for $\gamma>3$ with
$\gamma\not\in\N$. \\
a) Let $\gamma\in (1,3)\backslash\{2\}$. (Indeed the following arguments hold for any $\gamma\in (4\,n+1,4\,n+3)\backslash \N$ with
$n\in\N_0$.) Then $\tilde \alpha_0<0$ and for $z=i\,z_2$ with $z_2<0$, condition~(\ref{rootszabo}) simplifies to
\begin{equation*}
\begin{aligned}
     -\frac{z_2}{c_0}
        -\frac{1}{c_0}\,\mbox{Re}\left( \sqrt{ z_2^2  + 2\,\tilde\alpha_0\,c_0\,(-z_2)^{\gamma+1}}  \right)
   \leq C + N\,\log (1+|z_2|)\,.
\end{aligned}
\end{equation*}
Because $\gamma>1$ and $\tilde \alpha_0<0$, the term under the root is negative for sufficently large $-z_2$ and thus the real part of 
the root vanishes, which leads to the contradiction $\frac{|z_2|}{c_0}  \leq C + N\,\log (1+|z_2|)$. 
Therefore condition (B2) cannot be valid for any $\gamma\in(1,3)\backslash\{2\}$. \\
b) Now let $\gamma> 3$ and $z=z_1+i\,z_2$ with $z_1\in\R$ and $z_2<0$. We recall that 
$B(z)=\sqrt{1+ 2\,\tilde\alpha_0\,c_0\,(i\,z)^{\gamma-1}} - 1$. 
Let $\tilde z(r)=r\,e^{i\,\varphi}$ with $r>0$ and $\varphi:=\frac{\pi}{\gamma-1}-\frac{\pi}{2}$. Since $\gamma> 3$, we have $\varphi\in
(-\frac{\pi}{2},0)$ and thus $\tilde z_1>0$ and $\tilde z_2<0$. Since $(\gamma-1)\,(\varphi+\frac{\pi}{2})=\pi$ we have 
$\cos((\gamma-1)\,(\varphi+\frac{\pi}{2}))=-1$ and
$\sin((\gamma-1)\,(\varphi+\frac{\pi}{2}))=0$. This shows that 
\begin{equation}\label{helpszabo}
   1+ 2\,\tilde\alpha_0\,c_0\,(i\,\tilde z(r))^{\gamma-1}<0
  \qquad\mbox{for sufficiently large $r$.}
\end{equation}
Therefore $\mbox{Re}(B(\tilde z))=-1<0$, which together with $\tilde z_2<0$ implies $\tilde z_2\,\mbox{Re}(B(\tilde z))=-\tilde
z_2>0$. This shows that the first left hand side term of~(\ref{rootszabo}) is positive. 
Moreover,~(\ref{helpszabo}) implies that  $\mbox{Im}(B(\tilde
z))=\mbox{Im}(B(\tilde z)+1)>0$ for sufficiently large $-\tilde z_2$. From this together with $\tilde z_1>0$ 
we obtain $\tilde z_1\,\mbox{Im}(B(\tilde z))>0$, which shows that the second term on the left hand side of~(\ref{rootszabo}) 
is positive, too. If $r$ is sufficiently large then
$\mbox{Im}(B(\tilde z(r)))$ is of the order of $|r|^{\frac{\gamma-1}{2}}$ with $\gamma>3$ which cannot be bounded by $C + N\,\log (1+r)$.
Hence condition (B2) cannot be true for $\gamma>3$. In summary we have shown that the Green function of Szabo's equation cannot be
causal for $\gamma\in (1,\infty)\backslash\N$.\\
iii) Now we show that the front speed of $G_\alpha$ is $c_0$, if $\gamma\in (0,1)$.  
Since $G(\x,t+\frac{|\x|}{c_0})$ is causal, the front wave speed $v_\alpha$ satisfies $v_\alpha\leq c_0$. 
If $v_\alpha$ at $\x$ is smaller than
$c_0$, then $G(\x,t+\frac{|\x|}{c_0} +\delta\,|\x|)$ is causal for some $\delta> 0$.
This means that there exist constants $\epsilon>0$, $C>0$ and $N>0$ such that for all $z\in\R+i\,(-\infty,-\epsilon)$ 
(cf.~(\ref{coro:help1}) in Corollary~\ref{coro:powlaw1}):
\begin{equation*}
    -\mbox{Re}(iz\,(B(z)-\delta)) \leq C + N\,\log (1+|z|)\,,
\end{equation*}
since
\begin{equation*}
  \mbox{Re}(\alpha_*(-z)) + \mbox{Re}(i\,(-z)\,\delta) = -\mbox{Re}(iz\,(B(z)-\delta)) \,.
\end{equation*}
For the setting $z=i\,z_2$ with $z_2<0$, we obtain
\begin{equation*}
    \mbox{Re}\left( |z_2|\,\left[ \delta+1 -\sqrt{1 + 2\,\tilde\alpha_0\,c_0\, \left|\frac{1}{z_2}\right|^{1-\gamma}}\right] \right)
            \leq C + N\,\log (1+|z_2|)\,,
\end{equation*}
which cannot be true if $|z_2|$ is sufficiently large.  
Hence we conclude that $\delta=0$. This proves that the wave front speed is $c_0$ and concludes the
proof.
\end{proof}

\section{Causality analysis of the thermo-viscous wave equation}
\label{sec-caustv}

The operator of the \emph{thermo-viscous wave equation} (cf. e.g.~\cite{KiFrCoSa00})
$$
   P(D):=\left(\mbox{Id} +\tau_0\,\frac{\partial }{\partial t}\right)\,\nabla^2 -\frac{1}{c_0^2}\frac{\partial^2 }{\partial t^2}
$$
has order $3$ and \emph{principal part} $P_3(\mathbf{X}) =\tau_0\,t\,\sum_{j=1}^3 x_j$, where $\mathbf{X}:=(t,x_1,x_2,x_3)^T\in\R^4$. Since
$P_3(\mathbf{N})=0$ for $\mathbf{N}:=(1,0,0,0)^T$, the plane $\{ \mathbf{X}\in\R^4\,|\,\langle \mathbf{X},\mathbf{N}\rangle=0\}$ in $\R^4$ is
\emph{characteristic} with respect to P(D). According to Theorem 8.6.7 in~\cite{Ho03} the thermo-viscous wave equation with vanishing source
term has a solution $p_{tv}\in C^\infty(\R^4)$ such that $\mbox{supp}(p_{tv}) = \R^3\times\R^-_0$. This shows that the Green function of the
thermo-viscous wave equation is not uniquely determined. (The existence is guaranteed, since $P(\mathbf{X}):=\left(1
+\tau_0\,t\right)\,\sum_{j=1}^3 x_j^2 -\frac{t^2}{c_0^2}$ is not the zero polynomial.) Theorem~\ref{prop:thvis} below shows that the Green
function of the thermo-viscous wave equation cannot have a finite front speed and that a solution of the thermo-viscous wave equation depends
on its history. This explains why its Cauchy problem has no unique solution. For this theorem and Theorem~\ref{theo:mymodel} we need the
following lemma.

\begin{lemm}\label{lemm:T}
Let $\mathcal{S}_+=\{ f\in\mathcal{S}(\R) \,|\, \mbox{supp}(f)\subseteq \R^+_0 \}$ and $\gamma\in (1,2]$. The time-convolution operator $T_\gamma^\frac{1}{2}$ defined by the kernel
\begin{equation*}
\begin{aligned}
   k_{T_\gamma^{\frac{1}{2}}}(t)
      :=  \frac{1}{\sqrt{2\,\pi}}\,\mathcal{F}\left(\frac{1}{\sqrt{1+(-i\,\tau_0\,\omega)^{\gamma-1}}}\right)(t)
      %\qquad\mbox{ for $\gamma\in (1,2]$}
 \,
\end{aligned}
\end{equation*}
is an isomorphism of $\mathcal{S}_+$. Here the square root is understood as the primitive square
root. The inverse of $T_\gamma^{\frac{1}{2}}$ is the time-convolution operator $L_\gamma^\frac{1}{2}$ with the kernel
$$
k_{L_\gamma^{\frac{1}{2}}}(t) : =\frac{1}{\sqrt{2\,\pi}}\,\mathcal{F}\left(\sqrt{1+(-i\,\tau_0\,\omega)^{\gamma-1}}\right)(t) \,.
$$
Again the square root is understood as the primitive square root.
\end{lemm}

\begin{proof}
First we show that the kernel of $T_\gamma^{\frac{1}{2}}$ and $L_\gamma^{\frac{1}{2}}$ have supports in $[0,\infty)$. \\
%a) $\mbox{supp}\left( k_{T_\gamma^{\frac{1}{2}}}\right)\subseteq [0,\infty)$:
%b) $\mbox{supp}\left( k_{L_\gamma^{\frac{1}{2}}}\right)\subseteq [0,\infty)$:
Let $B(z):=\sqrt{1+(i\,\tau_0\,z)^{\gamma-1}}$ and $z\in \R+i\,\R^-$. Then
$$
\mbox{supp}\left(k_{T_\gamma^{\frac{1}{2}}}\right)\subseteq [0,\infty)   \quad\mbox{ and }\quad   \mbox{supp}\left(
k_{L_\gamma^{\frac{1}{2}}}\right)\subseteq [0,\infty)\,,
$$
if and only if $\hat u_1(z):=B(z)$ and $\hat u_2(z):=\frac{1}{B(z)}$ satisfy Theorem~\ref{th:hoer}. Since $\sqrt{1+(i\,\tau_0\,z)^{\gamma-1}}$
maps $\R+i\,\R^-$ analytically into $[1,\infty)+i\,\R$ and $1+(i\,\tau_0\,z)^{\gamma-1}$ cannot vanish on $\R+i\,\R^-$, $\hat u_1(z)$ and
$\hat u_2(z)$ map $\R+i\,\R^-$ analytically into $[1,\infty)+i\,\R$. Therefore property (A1) in Lemma~\ref{lemm:caus1} is satisfied. We
have for all $z\in \R+i\,\R^-$ with  $|z|>>1$:
\begin{equation*}
\begin{aligned}
         |B(z)| \leq C_1  \, |z|^{(\gamma-1)/2}
      \quad\mbox{ and }\quad
         \left|\frac{1}{B(z)}\right| \leq C_2  \, \left(\frac{1}{|z|}\right)^{(\gamma-1)/2} \qquad (\gamma\in (1,2])
\end{aligned}
\end{equation*}
for some constants $C_1,\,C_2>0$. Therefore property (A2) in Theorem~\ref{th:hoer} is satisfied for $\hat u_1(z)$ and $\hat u_2(z)$, which
proves that the kernels of $T_\gamma^{\frac{1}{2}}$ and $L_\gamma^{\frac{1}{2}}$ have support in $[0,\infty)$.

Since
$$
     T_\gamma^{\frac{1}{2}}\,f \in \mathcal{D}'_+ \qquad\mbox{ for every $f\in\mathcal{S}_+$}\,,
$$
and
$$
   \F^{-1}\left( T_\gamma^{\frac{1}{2}}\,f\right)(\omega) = \frac{\check f(\omega)}{\sqrt{1+(-i\,\tau_0\,\omega)^{\gamma-1}}} 
     \in \mathcal{S}  \qquad\mbox{ for every $f\in\mathcal{S}_+$}\,,
$$
the convolution operator $T_\gamma^{\frac{1}{2}}$ is well-defined on $\mathcal{S}_+$. 
Since $\mathcal{S}$ is invariant under multiplication by a polynomial, it follows analogously that $L_\gamma^{\frac{1}{2}}$ maps 
$\mathcal{S}_+$ into $\mathcal{S}_+$ and is well-defined. Since the Fourier transform is an isomorphism on $\mathcal{S}$ and
$$
      \F^{-1}\left(  L_\gamma^{\frac{1}{2}}\,T_\gamma^{\frac{1}{2}} \,f\right) = \check f 
= \F^{-1}\left(  T_\gamma^{\frac{1}{2}}\,L_\gamma^{\frac{1}{2}} \,f\right)\qquad\mbox{ for every $\check f\in\mathcal{S}$}\,,
$$
it follows that $L_\gamma^{\frac{1}{2}}:\mathcal{S}_+\to\mathcal{S}_+$ is the inverse of $T_\gamma^{\frac{1}{2}}:\X_+\to\mathcal{S}_+$ 
and $T_\gamma^{\frac{1}{2}}$ is an isomorphism of $\mathcal{S}_+$.
\end{proof}

With the help of the Laplace transform table in~\cite{Heu95} (cf. Appendix~2), the kernel of $T_2^{\frac{1}{2}}$ can be calculated 
as follows:
\begin{equation*}
\begin{aligned}
  k_{T_2^{\frac{1}{2}}}(t) &=
   \frac{1}{\sqrt{2\,\pi}}\,\mathcal{F}\left(\frac{1}{\sqrt{1-i\,\tau_0\,\omega}}\right)(t)
   = \sqrt{2\,\pi}\,\mathcal{L}^{-1}\left(\frac{1}{\sqrt{1+\tau_0\,p}}\right)(t) \\
   &= 2\,\sqrt{ \frac{\pi}{\tau_0\,t} }\,e^{-\frac{t}{\tau_0}}\,H(t)\,,
\end{aligned}
\end{equation*}
where $H(t)$ denotes the Heaviside function. 
This shows that for $\gamma=2$ the kernel of $T_2^{\frac{1}{2}}$ decreases exponentially. In the following we denote 
$T_\gamma^{\frac{1}{2}}\,T_\gamma^{\frac{1}{2}}$ by   $T_\gamma$.

\begin{prop}\label{prop:thvis}
Let $\tau_0,\,c_0>0$ be constants and $A(\omega):=1+\sqrt{1+(\tau_0\,\omega)^2}$ for any $\omega\in\R$. The Green function of the
thermo-viscous wave equation
\begin{equation}\label{thvwaveeq}
\begin{aligned}
    &\left(I+\tau_0\,\frac{\partial }{\partial t} \right)\,\nabla^2 G(\x,t)
       -\frac{1}{c_0^2}\frac{\partial^2 G(\x,t)}{\partial t^2}
            = -\delta(\x)\,\delta(t)
            \qquad  \R^3\times\R\,
\end{aligned}
\end{equation}
cannot have a finite front speed and is given by $G = T_2\,p_\alpha$, where $p_\alpha$ is defined as in~(\ref{patt}) with attenuation law
\begin{equation}\label{atthvis}
\begin{aligned}
   &\alpha(\omega)
      =\,\frac{\tau_0}{\sqrt{2\,A(\omega)}\,(A(\omega)-1)}\frac{\omega^2}{c_0}\,
\end{aligned}
\end{equation}
and phase speed
\begin{equation}\label{cthvis}
\begin{aligned}
 c(\omega)
              = \frac{\sqrt{2}\,(A(\omega)-1)}
          {\sqrt{A(\omega)}}\,c_0\,.
\end{aligned}
\end{equation}
\end{prop}

\begin{proof}
Applying the inverse Fourier transform to the thermo-viscous wave equation yields
\begin{equation}\label{help1}
\begin{aligned}
    &\nabla^2 \check G(\x,\omega)
       +k^2(\omega)\,\check G(\x,\omega)
            = -\frac{\delta(\x)}{\sqrt{2\,\pi}\,(1 -i\,\tau_0\,\omega)}  \qquad\mbox{ with}\\
   &k(\omega):= \frac{\pm\,\omega}{c_0\,\sqrt{1 - i\,\tau_0\,\omega}}\,.
\end{aligned}
\end{equation}
This problem has the solution $\check G(\x,\omega) = \frac{1}{\sqrt{2\,\pi}\,(1
-i\,\tau_0\,\omega)}\,\frac{e^{i\,k(\omega)\,|\x|}}{4\,\pi\,|\x|}$, where the square root of $1 - i\,\tau_0\,\omega$ is understood as the root
with positive real part. We assume that $G$ satisfies the causality requirement~(\ref{defg}) which in particular implies that
$\mbox{supp}(G(\x,\cdot))\in [0,\infty)$. Then the Green function can be written as follows
\begin{equation}\label{help2}
\begin{aligned}
      G(\x,t) =   T_2\,\left(\frac{1}{\sqrt{2\,\pi}}\int_\R
                   \frac{e^{-i\,(\omega\,t - k(\omega)\,|\x|)}}{4\,\pi\,|\x|} \,\d \omega\right)
         =: T_2\,p_\alpha (\x,t)\,,
%                     \qquad\mbox{ with $k(\omega)\in\C$}\,,
\end{aligned}
\end{equation}
where $T_2$ denotes the time-convolution operator in Lemma~\ref{lemm:T} for $\gamma=2$. The last identity is equivalent to
\begin{equation*}
\begin{aligned}
   \left( \mbox{Id} +\tau_0\,\frac{\partial}{\partial t} \right)   G(\x,t) =  p_\alpha (\x,t)\,,
\end{aligned}
\end{equation*}
which shows that $p_\alpha$ has finite front speed if and only if $G$ has finite front speed. 
Let $c_1>0$ be arbitrary but fixed. We prove a contradiction 
by showing that $p_\alpha(\x,t+\frac{|\x|}{c_1})$ cannot be a causal distribution for any $\x\in\R^3$. Comparing~(\ref{help2})
and~(\ref{patt}) shows that $\alpha(\omega)=\mbox{Im}(k(\omega))$ and $\frac{1}{c(\omega)}=\frac{\mbox{Re}(k(\omega)}{\omega}$. Since
$\alpha>0$ is required, the imaginary part of $k(\omega)$ must be positive and therefore we choose the positive sign for $k$, i.e.
\begin{equation}\label{kthv}
\begin{aligned}
    k(\omega) = \frac{\omega}{c_0}\frac{1}{A(\omega)-1}\left( \sqrt{\frac{A(\omega)}{2}}
                +i\frac{\tau_0\,\omega}{\sqrt{2\,A(\omega)}}\right)
\end{aligned}
\end{equation}
with $A(\omega):=1+\sqrt{1+(\tau_0\,\omega)^2}$. From this we get the attenuation law~(\ref{atthvis}) and the phase speed~(\ref{cthvis}). If
$p_\alpha(\x,t+\frac{|\x|}{c_1})$ is causal, then  property (B2) in Lemma~\ref{lemm:caus1} must be satisfied for $-i\,k(\omega) +
i\,\frac{\omega}{c_1}$, i.e. for each $\x\in\R^3$ there exist constants $\epsilon>0$, $C>0$ and $N>0$ such that for all $z\in
\R+i\,(-\infty,-\epsilon)$:
$$
    \mbox{Re}(i\,k(-z)) + \mbox{Re}\left( i\,\frac{z}{c_1}\right) \leq C + N\,\log (1+|z|)\,.
$$
For $z=i\,z_2$ with $z_2<0$ we get
$$
    \frac{1}{\sqrt{1-\tau_0\,z_2}}\,\frac{z_2}{c_0}
       - \frac{z_2}{c_1}
%     = -\mbox{Re}(\alpha_*(-z))
      \leq C + N\,\log (1+|z_2|)\,,
$$
which cannot be true for sufficiently large $-z_2$ and finite $c_1>0$. Hence property (B2) cannot be satisfied for 
$\alpha_*$ of the wave $p_\alpha$ if $c_1$ is finite. This contradiction proves that the front speed of $G$ cannot be finite and
concludes the proof.
\end{proof}

\section{Causal wave equations obeying attenuation power laws with $\gamma\in (1,2]$ for small frequencies}
\label{sec-gamgen}

As we have seen in Section~\ref{sec-causcheck}, the frequency power law for $\gamma\geq 1$ cannot hold. Now we present
admissible attenuation laws that permit approximate frequency power laws with $\gamma\in (1,2]$ for sufficiently small frequencies.\\

\noindent For given constants $\gamma\in (1,2]$, $0<c_1<\infty$ and $\tau_0\geq 0$ we define
\begin{equation}\label{alpha*mod}
\begin{aligned}
   \alpha_*(\omega) := \frac{-i\,\omega}{c_1\,\sqrt{1+(-i\,\tau_0\,\omega)^{\gamma-1}}}\,,
\end{aligned}
\end{equation}
where the square root is understood as the primitive square root. This implies for the attenuation law:
\begin{equation*}
\begin{aligned}
   \alpha(\omega)
      \approx \alpha_0\,|\tau_0\,\omega|^{\gamma}
\quad\mbox{ with } \quad
\alpha_0 = \frac{\sin(\frac{\pi}{2}\,(\gamma-1))}{\tau_0\,c_1}\,
\end{aligned}
\end{equation*}
for sufficiently small frequencies. Moreover, let $T^{\frac{1}{2}}_{\gamma}$ and $L^{\frac{1}{2}}_{\gamma}$ be defined as in
Lemma~\ref{lemm:T} and let the operators $D_*$ and $D_{*,|\x|}$ be defined as in Section~\ref{sec-genwaveeq}. Then
\begin{equation*}
\begin{aligned}
    D_* =
       \frac{1}{c_1}\,T^{\frac{1}{2}}_{\gamma}\,
       \frac{\partial}{\partial t}
\qquad\mbox{ and }\qquad
    D_{*,|\x|} = 0\,
\end{aligned}
\end{equation*}
and wave equation~(\ref{genwaveeq0}) (with $v_B$ replaced by $c_0$) reads as follows
\begin{equation}\label{waveeqtissue}
\begin{aligned}
%    \left(I+\tau_0^{\gamma-1}\,D_t^{\gamma-1} \right)\,\nabla^2 p_\alpha
%       -\left[\frac{c_0}{c_1}\,I+L^{\frac{1}{2}}_{\gamma}\right]^2\,
%          \frac{1}{c_0^2}\frac{\partial^2 p_\alpha}{\partial t^2}
%            = -f\,,
%
    \nabla^2 p_\alpha
       -\left[\mbox{Id} + \frac{1}{c_1}\,T^{\frac{1}{2}}_{\gamma}\right]^2\,
          \frac{1}{c_0^2}\frac{\partial^2 p_\alpha}{\partial t^2}
            = -f\,.
\end{aligned}
\end{equation}
For $\gamma=1$ we obtain the classical wave equation without damping and for $\gamma=2$, we obtain a modified thermo-viscous wave equation.
Since $\mbox{Re}(\alpha_*)$ is equal to~(\ref{atthvis}) for $\gamma=2$, the modified thermo-viscous wave equation obeys for $\gamma=2$ 
the same attenuation law as the thermo-viscous wave equation (if the source term $f$ is replaced by $L_2\,f$).

\begin{prop}\label{theo:mymodel}
The Green function of wave equation~(\ref{waveeqtissue}) has finite and constant wave front speed $c_0$.
\end{prop}

\begin{proof}
For the proof let $\x\in\R^3$ be arbitrary but fixed. The Green function $G(\x,t)$ of wave equation~(\ref{waveeqtissue}) has a front speed
$\leq c_0<\infty$, if $G_{tv}(\x,t+\frac{|\x|}{c_0})$ is causal, i.e. if $\alpha_*(\omega)$ defined by~(\ref{alpha*mod}) satisfies
Lemma~\ref{lemm:caus1}. We recall that the square root in the definition of $\alpha_*$ is understood as the primitive square root. According
to the proof of Lemma~\ref{lemm:T} $\frac{\alpha_*(-z)}{i\,z}$ satisfies property (B1) and thus $\alpha_*(-z)$ satisfies property (B1), too.

Property (B2) in Lemma~\ref{lemm:caus1} reads as follows: there exist constants $\epsilon>0$, $C>0$ and $N>0$ such that for all $z\in
\R+i\,(-\infty,-\epsilon)$:
\begin{equation}\label{condB1B2}
\begin{aligned}
    -\mbox{Re}\left(\alpha_*(-z)\right) =
    -\mbox{Re}\left(\frac{i\,z}{c_1\,\sqrt{1+(i\,\tau_0\,z)^{\gamma-1}}}\right)  \leq C + N\,\log (1+|z|)\,.
\end{aligned}
\end{equation}
Therefore (B2) is satisfied if $ -\mbox{Re}\left(\alpha_*(-z)\right)\leq 0$ for each $z\in \R+i\,\R^-$.
Let $B(z):=c_1\,\sqrt{1+(i\,\tau_0\,z)^{\gamma-1}}$ and $z=z_1+i\,z_2\in \R+i\,\R^-$ then
$$
    -\mbox{Re}\left(\alpha_*(-z)\right) =
     z_1\,\mbox{Im}\left(B^{-1}(z)\right) + z_2\,\mbox{Re}\left(B^{-1}(z)\right)\,.
$$
Since $\gamma\in (1,2]$ and $z_2<0$, it follows that $(i\,\tau_0\,z)^{\gamma-1}$ has positive real part and thus $B(z)$ has also positive real
part. The inversion of a complex number with positive real part yields a complex number with positive real part and thus
$\mbox{Re}\left(B^{-1}(z)\right)>0$  for any $z_2<0$. This proves that
\begin{equation}\label{B1}
\begin{aligned}
    z_2\,\mbox{Re}\left(B^{-1}(z)\right)<0      \quad\mbox{ for any $z\in \R+i\,\R^-$  }\,.
\end{aligned}
\end{equation}
For $z_1>0$ the imaginary part of $(i\,\tau_0\,z)^{\gamma-1}$ is positive and thus $\mbox{Im}\left(B(z)\right)>0$. The inversion of a complex
number with positive imaginary part yields a complex number with negative imaginary part and thus $\mbox{Im}\left(B^{-1}(z)\right)<0$.
%for any $z_1>0$.
Therefore we infer that $z_1\,\mbox{Im}\left(B^{-1}(z)\right)<0$ for any $z\in \R^++i\,\R^-$ with $z_1>0$. For $z_1<0$ the imaginary part of
$(i\,\tau_0\,z)^{\gamma-1}$ is negative and thus $\mbox{Im}\left(B(z)\right)<0$. Since the inversion of a complex number with negative
imaginary part yields a complex number with positive imaginary part we conclude that $\mbox{Im}\left(B^{-1}(z)\right)>0$.
% for any $z_1<0$.
Hence $z_1\,\mbox{Im}\left(B^{-1}(z)\right)<0$ for any $z\in \R^-+i\,\R^-$ with $z_1<0$. For $z_1=0$ we get
$z_1\,\mbox{Im}\left(B^{-1}(z)\right)=0$. In summary we have proven that
\begin{equation}\label{B2}
\begin{aligned}
   z_1\,\mbox{Im}\left(B^{-1}(z)\right)\leq 0
\quad\mbox{ for any $z\in \R+i\,\R^-$  }\,.
\end{aligned}
\end{equation}
(\ref{B1}) and~(\ref{B2}) imply that the left hand side of~(\ref{condB1B2}) is always non-positive and therefore ~(\ref{condB1B2}) is true.
This shows that $G(\x,t)$ has a front speed $\leq c_0<\infty$.

Now we show that the front speed of $G$ is equal to $c_0$. If the front speed is $v_\alpha(\x)<c_0$ for any $\x\in\R^3$,
then~(\ref{condB1B2}) must hold for
$$
     \alpha_*(-z) := \alpha_*(-z) + i\,\epsilon\,(-z)\,.
$$
For $z:=i\,z_2$ with sufficiently large  $-z_2$ we obtain
$$
     -\mbox{Re}\left(\alpha_*(-z)\right) = \frac{z_2}{c_1\,\sqrt{1+(-\tau_0\,z_2)^{\gamma-1}}}
       + \epsilon\,(-z_2)\,\,,
$$
which is positive and of the order $|z_2|$. This shows that condition (B2) can only be true if $\epsilon=0$. This concludes the proof.
\end{proof}

\begin{rema}
For $\gamma=2$ let $G_{c_0}$ denote the solution of wave equation~(\ref{waveeqtissue}) and let $G_{tv}$ denote the solution of the
thermo-viscous wave equation~(\ref{thvwaveeq}) with $c_0$ replaced by $c_1$. Then one can show that
$$
  \lim_{c_0\to\infty} L\,G_{c_0}(\x,t)
    = G_{tv}(\x,t)
  \qquad\mbox{ for each $\x\in\R^3$ and $t\in\R^+$}\,,
$$
which shows again that the front speed of $G_{tv}$ is infinite.
\end{rema}

\section{Acknowledgement}

This paper was partly supported by the "Frankreichschwerpunkt" of the University of Innsbruck. I would like to thank Xavier Bonnefond
and Pierre Mar\'echal for fruitful discussions.

%\bibliography{LitKoScBo08}

\bibliographystyle{plain}
\bibliography{LitKoScBo08}

\end{document}